
\input amstex.tex
\input amsppt.sty
\documentstyle{amsppt}
\let\e\varepsilon
\let\phi\varphi
\define\N{{\Bbb N}}

\define\A{{\Cal A}}

\def\ss{\subset}

\let\i\infty
\def\d{\delta}
\def\g{\gamma}

\def\a{\alpha}
\def\b{\beta}

\def\f{\frac}

\topmatter
\title On continuous   choice of  retractions  onto nonconvex  subsets
\endtitle

\author   Du\v san Repov\v s  and      Pavel V. Semenov   \endauthor
\leftheadtext{Du\v san Repov\v s  and Pavel V. Semenov}
\rightheadtext{ON CONTINUOUS CHOICE OF   RETRACTIONS}

\address
Faculty of Mathematics and Physics, and
Faculty of Education, 
University of
Ljubljana, P. O. B. 2964, Ljubljana, Slovenia 1001\endaddress
\email dusan.repovs\@guest.arnes.si \endemail

\address
Department of Mathematics, Moscow City Pedagogical University,
2-nd Selsko\-khozyast\-vennyi pr.\,4,\,Moscow,\,Russia 129226\endaddress
\email pavels\@orc.ru \endemail

\thanks
The first author was supported by the SRA grants
P1-0292-0101-04, J1-9643-0101 and  BI-RU/08-09-002.
The second author was supported by the
RFBR grant 08-01-00663. We thank the referee for comments and suggestions.
\endthanks

\abstract  For a Banach space $B$ and for a class $\A$ of its
bounded closed retracts, endowed with the Hausdorff metric, we
prove that  retractions on elements $A \in \A$ can be chosen to
depend continuously on $A$, whenever nonconvexity of each $A \in
\A$ is less than $\f{1}{2}$. The key geometric argument is that
the set of all uniform retractions onto an $\a-$paraconvex set
(in the
spirit of E. Michael) is $\frac{\a}{1-\a}-$paraconvex subset
in the space of continuous mappings of $B$ into itself. For a
Hilbert space $H$ the estimate $\frac{\a}{1-\a}$ can be improved
to $\frac{\a (1+\a^{2})}{1-\a^{2}}$ and the constant $\f{1}{2}$
can be reduced to the root of the equation $\a+ \a^{2}+a^{3}=1$.
\endabstract
\endtopmatter

\document

\head {\bf  0. Introduction}\endhead

The initial source of the present paper was two-fold. Probably it was
Bing \cite{12} who first asked whether there
exists a continuous function which selects a point from each arc
of the Euclidean plane. Hamstr\"{o}m and Dyer \cite{3}
observed that this problem reduces to the problem of
continuous choice of  retractions onto arcs. In fact, it suffices
to consider the images of a chosen point with respect to
continuously chosen retractions. A simple construction based, for
example, on the $\sin(\frac{1}{x})-$curve shows that in general
there are no continuously chosen retractions for the family of
arcs topologized by the Hausdorff metric. 

Therefore a stronger topology is needed
for an affirmative
answer. In fact, for any
homeomorphic compact subsets $A_{1}$ and $A_{2}$ of a metric space
$B$ one can consider the so-called $h-$metric $d_{h}(A_{1},
A_{2})$ defined by:
$$
d_{h}(A_{1}, A_{2})  =  \sup \{dist(x,h(x)): h \,\,\text{runs
over all homeomorphisms of $A_{1}$  onto $A_{2}$} \}
$$
and consider  the {\it completely regular} topology on the family
of all subarcs, generated by such a metric. With respect to this
topology, Pixley \cite{12} affirmatively resolved the problem of
continuous choice for retractions onto subarcs in an arbitrary
separable metric space.

By returning to the more standard Hausdorff topology in the
subspace $exp_{AR}(B)$ of all compact absolute retracts in $B$ one
can try to search for a degree of nonconvexity of such a retracts.
In the simplest situation, for convex exponent $exp_{conv}(B)$
consisting of all compact convex subsets of $B$, continuous choice
of retractions is a direct corollary of the following Michael
theorem \cite{8}:

\proclaim {Convex-Valued Selection Theorem} Any multivalued
mapping $F:X \to Y$ admits  a continuous singlevalued selection
$f:X \to Y,\,\,f(x) \in F(x)$,\,\,provided that:
\itemitem{(1)} $X$ is a paracompact space;
\itemitem{(2)} $Y$ is a Banach space;
\itemitem{(3)} $F$ is a lower semicontinuous (LSC) mapping;
\itemitem{(4)} For every $x \in X$,
$F(x)$ is a nonempty convex subset of $Y$; and
\itemitem{(5)} For every $x \in X$,
$F(x)$ is a closed subset of $Y$.
\endproclaim

In fact, let $X=exp_{conv}(B)$, let $Y$ be the  space $C_{b}(B,B)$ of
all continuous bounded mappings of $B$ into itself and
suppose that
$F:X \to Y$
associates to each $A \in X$ the nonempty set of all retractions
of $B$ onto $A$. Then all 
hypotheses $(1)-(5)$ can be verified
and the conclusion of the theorem gives the desired continuously
chosen retractions.

However, what  can one say about nonconvex absolute retracts? In
general, there exists an entire branch of mathematics devoted to
various generalizations and versions of the convexity. In our
opinion, even if one simply lists the titles of  "generalized
convexities" one will find as a minimum, nearly 20 different
notions. Among them are Menger's metric convexity \cite{7}, Levy's
abstract convexity \cite{5}, Michael's convex structures \cite{9},
Prodanov's algebraic convexity \cite{13}, M\"{a}gerl's paved
spaces \cite{6}, van de Vel's topological convexity \cite{21},
decomposable sets \cite{1}, Belyawski's simplicial convexity
\cite{2}, Horvath's structures \cite{4}, Saveliev's convexity
\cite{18}, and many others.

Typically, a creation of  "generalized convexities", is usually
related to an extraction of several principal properties of the
classical convexity which are used in one of the key mathematical
theorems or theories and, consequently deals with analysis and
generalization of these properties in maximally possible general
settings. Based on the ingenious idea of Michael who proposed in
\cite{10} the notion of a paraconvex set, the authors of
\cite{14-17,
19} systematically studied  another approach to weakening or
controlled omission of convexity on a set of principal theorems of
multivalued analysis and topology. Roughly speaking, to each
closed subset $P \ss B$ of a Banach space we have associated a
numerical function, say $\a_{P}:(0,+\i) \to [0,2)$, the so-called
function of nonconvexity of  $P$. The identity $\a_{P}\equiv 0$ is
equivalent to the convexity of $P$ and the more $\a_{P}$ differs
from zero the "less convex" is the set $P$.

Such classical results about multivalued mappings as the
 Michael
selection theorem, the
 Cellina approximation theorem, the Kakutani-Glicksberg fixed point theorem, the
 von Neumann - Sion minimax
theorem, etc. are valid with the replacement of the convexity
assumption for values $F(x),\,\,x \in X$ of a mapping $F$ by some
appropriate control of their functions of nonconvexity.

In comparison with usual ideas of "generalized convexities", we
never define in this approach, for example, a "generalized
segment" joining $x \in P$  and $y \in P$.
We look only for the
distances between points $z$ of the classical segment $[x,y]$  and
the set $P$ and look for the ratio of these distances and the size of
the segment. So the following natural question
immediately
arises:
Does paraconvexity of a set with respect to the classical
convexity structure coincide with convexity under some generalized
convexity structure?  Corollaries 2.5 and 2.6, based on continuous
choice of a retraction, in particular provide an affirmative
answer.

\head {\bf 1. Preliminaries}\endhead

Below we denote by $D(c,r)$ the open ball centered at the point $c$
with the radius $r$ and denote by $D_{r}$ an arbitrary open ball with
the radius $r$ in a metric space. So for a nonempty subset $P\ss
Y$ of a normed space $Y$, and for an open $r$-ball $D_r \ss Y$ we
define the relative precision of an approximation of $P$ by
elements of $D_r$ as follows:
$$ \delta (P,D_r) =
 \sup \left\{ \frac{dist(q,P)}{r}\, : \, q \in conv(P \cap D_r) \right\}.
$$

For a nonempty subset $P \ss Y$ of a normed space $Y$ the {\it
function $\alpha_{P}(\cdot)$ of nonconvexity} of $P$ associates to
each positive number $r$  the following nonnegative number
$$
\alpha_{P}(r) = \sup \lbrace \delta (P,D_r)\, \vert \, D_r
\,\,\text{is an open \,$r$-ball} \rbrace.
$$
Clearly, the identity $\alpha_{P}(\cdot) \equiv 0$ is equivalent
to the {\it convexity} of the closed set $P.$

\proclaim {Definition 1.1} For a nonnegative number $\a$ the
nonempty closed set $P$ is said to be $\a$-{\it paraconvex},
whenever $\a$ majorates the function $\a_{P}(\cdot)$ of
nonconvexity of the set, i.e
$$
dist(q,P) \leq \a \cdot r,\qquad \forall q \in conv(P \cap D_{r}).
$$
The nonempty closed set $P$ is said to be {\it paraconvex} if it
is $\a$-paraconvex for some $\a<1$.
\endproclaim

Recall, that a multivalued mapping $F: X \to Y$ between
topological spaces is called {\it lower semicontinuous} (LSC for
shortness) if for each open $U \subset Y$, its full preimage, i.e.
the set
$$
F^{-1}(U)\,=\,\{x \in X \vert\,\, F(x) \cap U \, \ne
\emptyset\,\}
$$
is open in $X$. Recall also that a singlevalued mapping $f:X \to
Y$ is called a {\it selection} (resp. an $\e$-{\it selection}) of
a multivalued mapping $F:X \to Y$ if $f(x) \in F(x)$\,(resp.\,
$dist(f(x),F(x)) < \e$\,), \, for all $x \in X$. Michael \cite{9}
proved the following selection theorem:

\proclaim {Paraconvex-Valued Selection Theorem} For each number
$0\leq \a <1$ any multivalued mapping $F:X \to Y$ admits  a
continuous singlevalued selection whenever:
\itemitem{(1)} $X$ is a paracompact space;
\itemitem{(2)} $Y$ is a Banach space;
\itemitem{(3)} $F$ is a lower semicontinuous (LSC) mapping; and
\itemitem{(4)} all values $F(x),\,\,x \in X$ are $\a-$paraconvex.
\endproclaim

As a corollary, every $\a-$paraconvex set, $\a <1$, is
contractible and moreover, it is an absolute extensor ($AE$) with
respect to the class of all paracompact spaces. Hence, it is an
absolute retract ($AR$). Moreover by \cite{17}, metric
$\e$-neighborhood of a paraconvex set in any uniformly convex
Banach space $Y$, is also a paraconvex set, and hence is also an
$AR$.

For each number $0\leq\a <1$ we denote by $exp_{\a}(B)$ the family
of all $\a-$paraconvex compact subsets and by  $bexp_{\a}(B)$ the
family of all $\a-$paraconvex bounded subsets of a Banach space
$B$ endowed with the Hausdorff metric. Recall that the Hausdorff
distance between two bounded sets is defined as the infimum of the
set of all $\e>0$ such that each of the sets is a subset of an
open $\e-$neighborhood of the other set.

For each retract $A \subset B$ we denote by $Retr(A)$ the set of
all continuous retractions of $B$ onto $A$. So the multivalued
mapping $Retr$ associates to each retract $A \subset B$ the set of
all retractions of $B$ onto $A$. For checking of the lower
semicontinuity of a mappings into the spaces of retractions and
for proving  paraconvexity of these
spaces we also need the notion
of an {\it uniform} retraction (in terminology of \cite{11}), or a
{\it regular} retraction (in terminology of \cite{20})). Recall
that a continuous retraction $R: B \to A$ is said to be {\it
uniform} (with respect to $A$) if
$$
\forall\,\, \e>0 \quad \exists \,\,\delta>0 \quad \forall \,\,x
\in B: \qquad dist(x, A)< \delta \Rightarrow dist(x, R(x)) < \e.
$$
We emphasize that a uniform retraction in general {\it is not} a
uniform mapping in the classical metric sense. Clearly, each
continuous retraction onto a {\it compact} subset is uniform with
respect to the set. So we denote by  $URetr(A)$ the set of all
continuous retractions of $B$ onto $A$ which are uniform with
respect to $A$.

\head {\bf 2.  The Banach space case}\endhead

\proclaim{Theorem 2.0} Let $0\leq \a < \f{1}{2}$ and $F: X \to
bexp_{\a}(B)$ be a continuous multivalued mapping of a paracompact
space $X$ into a Banach space $B$. Then there exists a continuous
singlevalued mapping $\goth{F}:X \to C_{b}(B,B)$ such that for
every $x \in X$ the mapping $\goth{F}_{x}:B \to B$ is a continuous
retraction of $B$ onto the value $F(x)$ of $F$.
\endproclaim
\demo{Sketch of proof of Theorem 2.0} Proposition 2.4 below is a
corollary of the Para\-convex\--valued selection theorem due to
Propositions 2.1-2.3 and the fact that $0 \leq\frac{\a}{1-\a}<1
\Leftrightarrow 0\leq \a <\f{1}{2}$. In turn, Theorem 2.0.
 follows directly from Proposition 2.4,
 it suffices to put
$\goth{F} = \goth{R}\circ F$.

\proclaim{Proposition 2.1} For every  $0\leq\a<1$ and for each
bounded $\a-$paraconvex subset $P$ the set $URetr(P)$ is a
nonempty closed subset of $C_{b}(B,B)$.
\endproclaim

\proclaim{Proposition 2.2} For every $0\leq\a<1$ and and for every
bounded $\a-$paraconvex subset $P \subset B$  the set $URetr(P)$
is an $\frac{\a}{1-\a}-$paraconvex subset of $C_{b}(B,B)$.
\endproclaim

\proclaim{Proposition 2.3} For every $0\leq\a<1$ the restriction
$URetr|_{bexp_{\a}(B)}$,\,\,$P\mapsto URetr(P)$ is lower
semicontinuous.
\endproclaim

\proclaim{Proposition 2.4} For every $0\leq\a<\f{1}{2}$ the
restriction $URetr|_{bexp_{\a}(B)}$,\,\, $P\mapsto URetr(P)$
admits a singlevalued continuous selection, say
$$\goth{R}:
bexp_{\a}(B) \to C_{b}(B,B),\,\,\goth{R}_{P} \in URetr(P).$$
\endproclaim
\enddemo

\demo{Proof of Proposition 2.1} Clearly for each bounded closed
retract $A$ the sets $Retr(A)$ and $URetr(A)$ are closed in the
Banach space $C_{b}(B,B)$. To obtain the nonemptiness of $Retr(P)$
for the $\a-$paraconvex set $P$ it suffices to apply the
Paraconvex-valued selection theorem to the mapping $F:B \to B$
defined by setting $F(x)=P$ for $x \in B \setminus P$ and
$F(x)=\{x\}$ for $x \in P$. To construct a uniform retraction $R:
B \to P$ one must study more in detail the idea of the proof of
the Paraconvex-valued selection theorem.

Let us denote by $d(x)$ the distance between a point $x \in B$ and
a fixed $\a-$para\-convex subset $P \subset B$. For every $x \in
B\setminus P$ first consider the intersection of the set $P$ with
the open ball $D(x, 2d(x))$. Next, take the convex hull $conv\{P
\cap D(x, 2d(x))\}$ and finally, define the convex-valued mapping
$H_{1}: B\setminus P \to B$ by setting
$$
H_{1}(x) = \overline{conv}\{P \cap D(x, 2d(x))\}.
$$
This mapping is a LSC mapping defined on the paracompact domain
$B\setminus P$ with nonempty closed convex values in Banach space.
So the  Convex-valued selection theorem guarantees the existence
of a continuous singlevalued selection, say $h_{1}: B\setminus P
\to B,\,\,h_{1}(x) \in H_{1}(x)$.

For an arbitrary $\a <\b < 1$ the $\a-$paraconvexity of $P$
implies the inequalities
$$
dist(h_{1}(x), P) <  \b \cdot2d(x),\qquad dist(x,h_{1}(x)) \leq
2d(x),\qquad x \in B\setminus P.
$$

Similarly, define the convex-valued and closed-valued LSC mapping
$H_{2}: B\setminus P \to B$ by setting
$$
H_{2}(x) = \overline{conv}\{P \cap D(h_{1}(x), \b \cdot
2d(x))\},\qquad x \in B\setminus P.
$$
For its continuous singlevalued selection $h_{2}: B\setminus P \to
B,\,\,h_{2}(x) \in H_{2}(x)$ we see that for every
$x \in B\setminus P$,
$$
dist(h_{2}(x), P) \leq\a \cdot \b \cdot 2d(x) <
\b^{2}\cdot2d(x),
\quad
dist(h_{2}(x), h_{1}(x)) \leq \b
\cdot2d(x),$$
once again due to the $\a-$paraconvexity of $P$.

One can inductively construct a sequence
$\{h_{n}\}_{n=1}^{\infty}$ of continuous singlevalued mappings
$h_{n}: B\setminus P \to B$ such that for every
$x \in B\setminus P$,
$$
dist(h_{n+1}(x), P) <  \b^{n+1}\cdot2d(x),\quad dist(h_{n+1}(x),
h_{n}(x)) \leq \b^{n}\cdot2d(x).
$$

So the sequence $\{h_{n}\}_{n=1}^{\infty}$  is locally uniformly
convergent and hence its pointwise limit $h(x)= lim_{n \to
\infty}\,h_{n}(x)$ is well-defined and continuous. Moreover, $h(x)
\in P,\,\,x \in B\setminus P$, due to the closedness of $P$ and
convergency of $\{h_{n}\}_{n=1}^{\infty}$.

Hence the mapping $R : B \to P$ defined by $R(x)=h(x),\,\,x \in
B\setminus P$ and $R(x)=x,\,\, x \in P$ is a retraction of $B$
onto $P$ which is continuous over the set $B\setminus P$ by
construction.

To finish the proof we
estimate that for every $x \in
B\setminus P$:
$$dist(x,h(x))
\leq dist(x,h_{1}(x))+
\sum_{n=1}^{\infty}\,dist(h_{n}(x),h_{n+1}(x))\leq$$
$$\leq 2d(x)(1+\b+\b^{2}+\b ^{3}+...) = C\cdot d(x),
$$
for the constant $C=\frac{2}{1-\b}$. So  for $x_{0} \in P$ and for
$x \in B\setminus P$ we have
$$
dist(R(x_{0}),R(x))=dist(x_{0},h(x))\leq
dist(x_{0},x)+dist(x,h(x))\leq
$$
$$
\leq dist(x_{0},x)+C\cdot d(x) \leq (1+C)dist(x_{0},x).
$$

The continuity of the retraction $R: B \to P$ over the closed
subset $P \subset B$ and its uniformity clearly follow from the
last inequality. \qed
\enddemo

\medskip

\demo{Proof of Proposition 2.2} Pick an open ball $D(h,r)$ with
the radius $r$ in the space $C_{b}(B, B)$ centered at the mapping
$h \in C_{b}(B, B)$ which intersects with the closed set
$URetr(P)$. Let $R_{1}, R_{2},...,R_{n}$  be elements of the
intersection $D(h,r) \cap URetr(P)$ and let $Q$ be a convex
combination of  $R_{1}, R_{2},...,R_{n}$. We want to estimate the
distance between $Q$ and $URetr(P)$.

Pick a point $x \in B$. Passing from the mappings $h,Q,R_{1},
R_{2},...,R_{n} \in C_{b}(B, B)$ to their values at $x$ we find
the open ball $D(h(x),r)$ with the radius $r$ in the space $B$
centered at $h(x) \in B$, the finite set $\{R_{1}(x),
R_{2}(x),...,R_{n}(x)\}$ of elements from the intersection
$D(h(x),r)\bigcap P$ and the point $Q(x) \in conv(D(h(x),r)\cap
P)$. Having all fixed continuous mappings $h,Q,R_{1},
R_{2},...,R_{n} \in C_{b}(B, B)$ we see that all points
$h(x),Q(x),R_{1}(x), R_{2}(x),...,R_{n}(x) \in B$ continuously
depend on $x \in B$.

Let $r(x)$ be the Chebyshev radius of the compact convex
finite-dimensional set
$$\Delta(x)= conv\{R_{1}(x),...,R_{n}(x)\},$$ i.e. the infimum (in fact, the minimum),
of the set of radii of all closed balls containing $\Delta(x)$.
Clearly $r(x) < r,\,\,x \in X$. Moreover $r(x)$ continuously
depends on $x$ and for any positive $\g >0$ the entire set
$\Delta(x)$ lies in the open ball $D(C(x),r(x)+\g)$ for some
suitable point $C(x) \in \Delta(x)$.

Henceforth, the $\a-$paraconvexity of $P$ implies that for an
arbitrary $\a<\b$ the inequality 
$$dist(Q(x),P)<\b \cdot
\varrho(x), \ \ \ \ \varrho(x)=r(x)+\gamma$$ 
holds. So, the multivalued mapping
$$
F_{1}(x) = \overline{conv}\{P \cap D(Q(x), \b \cdot \varrho(x))
 \}.
$$
is a LSC mapping with nonempty convex and closed values. Note that
for each $x \in P$ all points
$R_{1}(x),R_{2}(x),...,R_{n}(x),Q(x)$ coincide with $x$ because
all $R_{1},...,R_{n}$ are retractions onto $P$. So the identity
mapping $id|_{P}$ is a continuous selection of $F_{1}|_{P}$.
Therefore the mapping $G_{1}$ which is identity on $P \subset B$
and otherwise coincides with $F_{1}$  admits a continuous
singlevalued selection, say $Q_{1}: B \to B,\,\,Q_{1}(x) \in
G_{1}(x)$.The $\a-$paraconvexity of $P$ and the construction imply
that
$$
dist(Q_{1}(x),P)<\b^{2} \cdot \varrho(x),\qquad
dist(Q_{1}(x),Q(x))\leq \b \cdot \varrho(x),\qquad
Q_{1}|_{P}=id|_{P}.
$$

Similarly, the multivalued mapping defined by setting
$$
F_{2}(x)= \overline{conv}\{P \cap D(Q_{1}(x), \b^{2}\cdot
\varrho(x))\}
$$
admits a continuous singlevalued selection, say $Q_{2}: B \to B$
such that
$$
dist(Q_{2}(x),P)<\b^{3} \cdot \varrho(x),\qquad
dist(Q_{2}(x),Q_{1}(x))\leq \b^{2} \cdot \varrho(x),\qquad
Q_{2}|_{P}=id|_{P}.
$$

Inductively we obtain a sequence $\{Q_{n}\}_{n=1}^{\infty}$ of
continuous singlevalued mappings $Q_{n}: B \to B$ with the
properties that $Q_{n}|_{P}=id|_{P}$ and
$$
dist(Q_{n+1}(x), P) <  \b^{n+2}\cdot \varrho(x), \qquad
dist(Q_{n+1}(x), Q_{n}(x)) \leq \b^{n+1}\cdot \varrho(x).
$$

Clearly the pointwise limit $R$ of the sequence
$\{Q_{n}\}_{n=1}^{\infty}$ is continuous retraction of $B$ onto
$P$ and, moreover,
$$
dist(Q(x),R(x)) \leq dist(Q(x),
Q_{1}(x))+\sum_{n=1}^{\infty}\,dist(Q_{n}(x),Q_{n+1}(x))\leq
$$
$$
\leq \b \cdot(1+\b+\b^{2}+\b^{3}+...)\cdot \varrho(x)
=\frac{\b}{1-\b}\cdot \varrho(x).
$$
Hence,
$$dist(Q, Retr(P)) \leq \frac{\b}{1-\b}\cdot \varrho(x) =
\frac{\b}{1-\b}\cdot (r(x)+\g) < \frac{\b}{1-\b}\cdot (r+\g).
$$
Passing to $\b \to \a+0$ and to $\g \to 0+0$ we conclude $dist(Q,
Retr(P))\leq \frac{\a}{1-\a} \cdot r$. To finish the proof we must
check that the retractions $R(x)=\lim_{n \to \i}Q_{n}(x),\,\,x\in
X$ onto $P$ constructed above are uniform with respect to $P$. To
this end, using uniformity of all retractions $R_{1},...,R_{n}$,
for an arbitrary $\e>0$ choose $\d>0$ such that
$$
 dist(x, P)< \delta \Rightarrow dist(x, R_{i}(x)) < \e.
$$
In particular, for every  point $x$ with $dist(x, P)< \delta$ all
values $R_{1}(x),...,R_{n}(x), Q(x)$ are in the open ball $D(x,
\e)$. Hence $r(x) < \e$ and this is why
$$
dist(x,R(x)) \leq dist(x,Q(x)) + dist(Q(x),R(x)) < \e +
\frac{\b}{1-\b}\cdot \varrho(x) < \frac{1}{1-\b}\cdot (\e+\g)
$$
Therefore $R \in URetr(P)$ and $dist(Q, URetr(P)) \leq
\frac{\a}{1-\a}\cdot r$. So $URetr(P)$ is
$\frac{\a}{1-\a}-$paraconvex. \qed  \enddemo

\medskip

\demo{Proof of Proposition 2.3} Pick $P \in bexp_{\a}(B)$, an
uniform retraction $R \in URetr(P)$ and a number $\e > 0$. So  let
$\delta > 0$ be such that $\delta < (1-\a)\cdot \e$ and
$$
dist(x, P)< \delta \quad \Rightarrow \quad dist(x, R(x)) <
(1-\a)\cdot\e.
$$

Consider any $P'\in bexp_{\a}(B)$ which is $\delta-$close to $P$
with respect to the Hausdorff distance. We must find a uniform
retraction $R' \in URetr(P')$ such that $dist(R,R')< \e.$

The multivalued mapping $F':B \to B$ such that $F'(x)=\{x\},\,\,x
\in P'$ and $F'(x)=P'$ otherwise is a LSC mapping with
$\a-$paraconvex values. Any selection of $F'$ is a retraction onto
$P'$. So let us check that $R$ is almost selection of $F'$ and
hence, is almost a retraction onto $P'$.

For every $x \in B \setminus P'$ we have
$$dist(R(x),F'(x))=dist(R(x),P')< \delta < (1-\a)\e$$
because $R(x) \in P$ and the set $P$ lies in the
$\delta-$neighborhood of the set $P'$. If $x \in P'$ then
$$dist(R(x),F'(x))=dist(R(x),x)<(1-\a)\e$$
because the set $P'$ lies in the $\delta-$neighborhood of the set
$P$ and due to the choice of the number $\delta$. Hence, the
retraction $R$ of $B$ onto the set $P$ is a continuous
singlevalued $\e'-$selection of the mapping $F'$,\,\,
$\e'=(1-\a)\e$.

Following the proofs of Propositions 2.1 and 2.2 we can improve
the $\e'-$selection $R$ of $F'$ to a selection $R'$ of $F'$ such
that $dist(R,R')<\frac{\e'}{1-\a}=\e$. So $R'$ is a continuous
retraction onto $P'$ and the checking of uniformity of $R'$ can be
performed by repeating the arguments on Chebyshev radii 
from the proof of Proposition 2.2. \qed
\enddemo
\medskip

Observe the proof of Theorems 2.0 for the case of compact
paraconvex sets is much more easier, because for any compact
retract $A \ss B$ each continuous retraction $B \rightarrow A$
automatically will be uniform with respect to $A$. So, one can
uses directly $Retr(A)$ instead of $URetr(A)$.

\proclaim{Corollary 2.5}  Under the assumptions of Theorem 2.0 if
in addition all values $F(x),\,\,x \in X$, are pairwise disjoint
then the metric subspace $Y= \bigcup_{x \in X}F(x) \ss B$ admits a
convex metric structure $\sigma$ (in the sense of
\cite{9})\,\,such that each value $F(x)$ is convex with respect to
$\sigma$.
\endproclaim
\demo{Proof} By Theorem 2.0, let $R(x): B \to F(x),\,\,x \in X$,
be a continuous family of uniform continuous retractions onto the
values $F(x)$. One can define a convex metric structure $\sigma$
on $Y= \bigcup_{x \in X}F(x)$ by setting that $\sigma-$convex
combinations are defined only for finite subsets
$\{y_{1},y_{2},...,y_{n}\}$ which are entirely displaced in a
value $F(x)$ and
$$\sigma-conv_{F(x)}\{y_{1},y_{2},...,y_{n}\}=
R(x)(conv_{B}\{y_{1},y_{2},...,y_{n}\}).\qed 
$$
\enddemo

\proclaim{Corollary 2.6} Let $f: Y \to X$ be a continuous
singlevalued surjection and let all point-inverses
$f^{-1}(x),\,\,x \in X$ are $\a-$paraconvex subcompacta of $Y$
with $\a < \f{1}{2}$. Then $Y$ admits a convexity metric structure
such that each point-inverse is convex with respect to this
structure.
\endproclaim
\medskip

\head {\bf 3. The Hilbert space case}\endhead

Hilbert spaces have a many of advantages inside the class of all
Banach spaces. In this chapter we demonstrate such a advantage
related to paraconvexity. Briefly we prove the estimate
$\frac{\a}{1-\a}$ for paraconvexity of the set $Retr(P)$ onto
$\a-$paraconvex set $P$ can be improved with
$$
\frac{\a (1+\a^{2})}{1-\a^{2}} = \frac{\a}{1-\a} \cdot
\frac{1+\a^{2}}{1+\a} < \frac{\a}{1-\a}. $$

Hence in  Theorem 2.0 one can substitute the root of the equation
$\a+\a^{2}+\a^{3}=1$ instead of $\f{1}{2}$.  In fact, a
generalization of such type can be performed  for any uniformly
convex Banach spaces but  for Hilbert space
the proofs differ only in
technical details.

\proclaim{Theorem 3.0}  Let $H$ be a Hilbert space and $F: X \to
bexp_{\a}(H)$ be a continuous mapping of a paracompact space $X$,
where $\a+\a^{2}+\a^{3}<1$. Then there exists a continuous
singlevalued mapping $\goth{F}:X \to C_{b}(H,H)$ such that for
every $x \in X$ the mapping $\goth{F}_{x}:H \to H$ is a continuous
retraction of $H$ onto the value $F(x)$ of $F$.
\endproclaim

So  we repeat the original definition of $\a-$paraconvexity of $P$
but with the appropriate estimate for distances between points of
simplices and points of $P$ {\it inside}  open balls.

\proclaim {Definition 3.1} Let $0\leq \a < 1$. A nonempty closed
subset $P \subset B$ of a Banach space $B$ is said to be {\rm
strongly} $\a$-{\rm paraconvex} if for every open ball $D \subset
B$ with radius $r$ and for every $q \in conv(P \cap D)$ the
distance $dist(q,P \cap D)$ is less than or equal to $\a \cdot r$.
\endproclaim

Clearly, strong $\a-$para\-convexity of a set implies its
$\a-$para\-convexity. In a Hilbert space the converse is almost
true: for some
$1>\b>\a$,  $\a-$para\-convexity implies strong
$\b-$\-para\-convexity for some suitable $\b$.

\proclaim{Proposition 3.2} Any $\a-$paraconvex subset $P$ of a
Hilbert space is its strong $\varphi(\a)-$paraconvex subset, where
$\varphi(\a) = \sqrt{1-(1-\a)^{2}}=\sqrt{2\a - \a^{2}}$.
\endproclaim

Proposition 3.2 is an immediate corollary of the following purely
geometrical lemma:

\proclaim{Lemma 3.3} Let $D=D_{r}$ be an open ball with the radius
$r$ in the Hilbert space $H$. Let $z$ be a point of the convex
hull $conv(P \cap D)$ of the intersection $D$ with a set $P$ and
let $dist(z,P) \leq \a \cdot r$. Then $dist(z,P \cap D) \leq
\varphi(\a) \cdot r$
\endproclaim
\demo{Proof of Lemma 3.3}  Pick an arbitrary $\a<\g<1$ and let $c$
be the center of the open ball $D=D(c,r)$.

If $dist(c,z)\leq (1-\g)\cdot r$ then the whole open ball $D(z, \g
\cdot r)$ lies inside of $D$. Hence, a point $p \in P$ which is
$(\g \cdot r)-$close to $z$ automatically will be in $D$. So
$$
dist(z, P \cap D) \leq dist(z,p) < \g \cdot r \leq \varphi(\g)
\cdot r.
$$

Let us look for the opposite case when $z$ is "close" to the
boundary of the ball $D$, i.e. when $(1-\g)\cdot r < dist(c,z) < r
$. Draw the hyperplane $\Pi$  supporting to the ball
$D(c,\,dist(c,z))$ at the point $z$. Such the hyperplane $\Pi$
divides the ball $D$ into two open convex parts: the center $c$
belongs to the "larger" part $D_{+}$ whereas the point $z$ belongs
to the the boundary of "smaller" part $D_{-}$. Clearly,
$Clos(D_{-})$ contains a point $p \in P$ (if, to the contrary, $P
\cap D$ is subset of $D_{+}$ then $z \in conv(P \cap D) \subset
D_{+}$). Hence, the distance $dist(z, P \cap D )$ majorates by
$$
dist(z,p) \leq \max \{dist(z,u): \,\,u \in Clos(D_{-})\} = \varphi
\left(\frac{dist(c,z)}{r}\right) \cdot r< \varphi(\g) \cdot r.
$$

So  in both cases $dist(z, P \cap D) \leq \varphi(\g) \cdot r$ and
passing to $\g \to \a+0$ we see that $dist(z, P \cap D) \leq
\varphi(\a) \cdot r$ \qed
\enddemo

Recall that for a multivalued mapping $F: X \to Y$ and for a
numerical function $\e: X \to (0,+\i)$ a singlevalued mapping $f:X
\to Y$ is said to be an $\e-$selection of $F$ if $dist(f(x),F(x))<
\e(x),\,\,x \in X$.

\proclaim{Proposition 3.4} Let $0\leq \a <1$ and let $F: X \to H$
be an $\a-$paraconvex valued LSC mapping from a paracompact domain
into a Hilbert space. Then
\itemitem{(1)} for each constant $C > \frac{1+\a^{2}}{1-\a^{2}}$,
for every continuous function $\e : X \to (0,+\i)$ and for every
continuous $\e-$selection $f_{\e}:X \to H$ of the mapping $F$
there exists a continuous selection $f:X \to H$ of $F$ such that
$$
 dist(f_{\e}(x), f(x)) < C \cdot \e(x), \qquad x \in X;
$$
\itemitem{(2)} $F$ admits a continuous selection $f$.
\endproclaim
\demo{Proof} Clearly $(1)$ implies $(2)$: the mapping $x \mapsto
[1+dist(0,F(x)),+\i)$,\,\,$x \in X$, is a LSC mapping with
nonempty closed and convex  values and therefore admits a
continuous selection, say $\e: X \to (0,+\i)$.
Therefore $f_\e \equiv 0$ is an $\e$-selection of $F$.

To prove $(1)$ let $\varphi(t) = \sqrt{2t - t^{2}},\,\, 0<t<1$,
choose any $\g \in (\a,1)$ and denote by $D(x)=D(f_{\e}(x),
\e(x))$. As above, the multivalued mapping
$$
F_{1}(x) = \overline{conv}\{F(x) \cap D(x) \},\qquad x \in X
$$
admits a singlevalued continuous selection, say $f_{1}:X \to H$.

For each $x \in X$ the point $f_{1}(x)$ belongs to the convex hull
$conv\{F(x) \cap D(x)\}$ and $dist(f_{1}(x),F(x)) \leq \a \cdot
\e(x)$ due to the $\a-$paraconvexity of the value $F(x)$. Lemma
3.3 implies that
$$
dist(f_{1}(x),F(x) \cap D(x)) \leq \varphi(\a) \cdot \e(x) <
\varphi(\g) \cdot \e(x).
$$

Therefore, the multivalued mapping $F_{2}:X \to H$ defined by
$$
F_{2}(x) = \overline{conv}\{F(x) \cap D(x) \cap D(f_{1}(x),
\varphi(\g) \cdot \e(x))\},\qquad x \in X
$$
is a LSC mapping with nonempty closed and convex values. Hence
there exists a selection of $F_{2}$, say $f_{2}:X \to H$.

For each $x \in X$ the point $f_{2}(x)$ belongs to the convex hull
$\overline{conv}\{F(x) \cap D(x)\}$ and $dist(f_{2}(x),F(x)) \leq
\a \cdot \varphi(\g) \cdot \e(x)$ due to the $\a-$paraconvexity of
the value $F(x)$ and because $f_{2}(x) \in \overline{conv}\{F(x)
\cap D(f_{1}(x), \varphi(\g) \cdot \e(x))\}.$

Lemma 3.3 implies that
$$
dist(f_{2}(x), F(x)) \leq \varphi(\a \cdot \varphi(\g)) \cdot
\e(x) < \varphi(\g \cdot \varphi(\g)) \cdot \e(x),\qquad x \in X.
$$

Put
$$
F_{3}(x) = \overline{conv}\{F(x) \cap D(f_{\e}(x), \e(x)) \cap
D(f_{2}(x), \varphi(\g \cdot \varphi(\g)) \cdot \e(x))\},\qquad x
\in X
$$
and so on. Hence we have constructed a sequence $\{f_{n}:X \to
H\}^{\i}_{n=1}$ of continuous singlevalued mappings such that
$$
dist(f_{\e}(x), f_{n}(x)) \leq \e(x),\qquad dist(f_{n}(x),F(x)) <
\g_{n} \cdot \e(x)
$$
where $\g_{1}=\g$ and $\g_{n+1}=\g \cdot \varphi(\g_{n}).$

The sequence $\{\g_{n}\}$ is monotone, decreasing and converges to
the (positive!) root of equation $t=\g \cdot \varphi(t)$, i.e. to
the number $t=\frac{2\g^{2}}{1+\g^{2}} >
\frac{2\a^{2}}{1+\a^{2}}$. Therefore we can choose numbers $N \in
\N$ and $\lambda$ such that
$$
1> 1- \frac{1}{C} > \lambda> \g_{N}>\frac{2\g^{2}}{1+\g^{2}} >
\frac{2\a^{2}}{1+\a^{2}}.
$$
Hence, the mapping $g_{1}=f_{N}$ is a $(\lambda \cdot
\e)-$selection of $F$ and $$dist(f_{\e}(x)), g_{1}(x))\leq
\e(x).$$

Starting with $g_{1}$ one can find $\lambda^{2} \cdot
\e-$selection $g_{2}$ of $F$ such that
$$dist(g_{1}(x), g_{2}(x))\leq \lambda \cdot \e(x).$$
Continuation of this construction produces a continuous selection
$f=lim_{n \to \i} g_{n}$ of $F$ such that
$$
dist(f_{\e}(x), f(x)) \leq \e(x) \cdot
(1+\lambda+\lambda^{2}+...)= \frac{1}{1-\lambda} \cdot \e(x) < C
\cdot \e(x), \qquad x \in X.
$$
\qed
\enddemo
\medskip

Proposition 3.4 implies the following analog of Proposition 2.2:

\proclaim{Corollary 3.5} For every $0\leq\a<1$ and for every
bounded $\a-$paraconvex subset $P \subset H$  the set $Retr(P)$ is
an $\frac{\a (1+\a^{2})}{1-\a^{2}}-$paraconvex subset of
$C_{b}(H,H)$.
\endproclaim
\medskip

\demo{Proof of Theorem 3.0} It suffices to repeat the proof of
Theorem 2.0 but we use Corollary 3.5 instead of Proposition 2.2.
\qed
\enddemo

\head {\bf 4. Concluding remarks}\endhead

Roughly speaking, we have proved that $\a-$paraconvexity of a
set implies $\beta-$para\-convexity of a set of all retractions onto
this set with $\beta=\beta(\a)=\frac{\a}{1-\a}$. Such an estimate
for $\beta=\beta(\a)$ naturally appears as a result of the
usual
geometric progression procedure. However, it is unclear to us
whether the constant $\frac{\a}{1-\a}$ is the best possible?

Some examples in the Euclidean plane show that in some particular
cases (for some curves) the constant $\beta=\beta(\a)$ admits more
precise estimates. Unfortunately, calculations in these examples
are based on geometric properties of concrete $\a-$ paraconvex
curves which are in general false for arbitrary $\a-$paraconvex
subsets of the plane.

Hence the question about continuous choice of retractions onto
bounded $\a$-paraconvex sets with $\f{1}{2}\leq \a <1$ remains
open. Even the case of subsets of the Euclidean plane presents
an evident interest. The main obstructions for a various
counterexamples are related to problems of constructing
retractions with some prescribed constraints.


\Refs

\widestnumber\key{1000000}

\ref \key{1}
\by H. A. Antosiewicz and A. Cellina
\paper Continuous selections and differential relations
\jour J. Diff. Eq.
\vol {\bf 19 }
\yr 1975
\pages 386--398
\endref

\ref \key{2}
\by R. Bielawski
\paper Simplicial convexity and its applications
\jour J. Math. Anal. Appl.
\vol {\bf 106 }
\yr 1986
\pages 155--171
\endref

\ref \key{3}
\by M.-E. Hamstr\"{o}m and E. Dyer
\paper Completely regular mappings
\jour Fund. Math.
\vol {\bf 45 }
\yr 1958
\pages 103--118
\endref

\ref \key{4}
\by C. D. Horvath
\paper Contractibility and generalized convexity
\jour J. Math. Anal. Appl.
\vol {\bf 156}
\yr 1991
\pages 341--357
\endref

\ref \key{5} 
\by F. W. Levy 
\paper On Helly's theorem and the axioms of convexity 
\jour J. Indian Math. Soc. 
\vol{\bf 15 } 
\yr 1951 
\pages 65--76
\endref

\ref \key{6}
\by G. M\"{a}gerl
\paper An unified approach to measurable and continuous selections
\jour Trans. Amer. Math. Soc.
\vol {\bf 245 }
\yr 1978
\pages 443--452
\endref

\ref \key{7}
\by K. Menger
\paper Untersuchungen uber allgemeine Metrik. I, II, III
\jour Math. Ann.
\vol {\bf 100 }
\yr 1928
\pages 75--163
\endref

\ref \key{8}
\by E. Michael
\paper Continuous selections. I
\jour Ann. of Math.
\vol{\bf (2) 63 }
\yr 1956
\pages 361--382
\endref

\ref \key{9}
\by  E. Michael 
\paper Convex structures and continuous selections
\jour Canadian J. Math.
\vol {\bf 11}
\yr 1959
\pages 556--575
\endref

\ref \key{10}
\by  E. Michael
\paper Paraconvex sets
\jour Math. Scand.
\vol {\bf 7 }
\yr 1959
\pages 372--376
\endref

\ref \key{11}
\by  E. Michael
\paper Uniform AR's and ANR's
\jour Compositio Math.
\vol {\bf 39}
\yr 1979 \pages 129--139
\endref

\ref \key{12}
\by C. Pixley
\paper Continuously choosing a retraction of a separable metric space onto each its arcs
\jour Illinois J.  Math.
\vol {\bf 20}
\yr 1976
\pages 22--29
\endref

\ref \key{13} \by I. Prodanov 
\paper A generalization of some
separability theorems 
\jour C. R. Acad. Bulgare. Sci. 
\vol {\bf 17}
\yr 1964 
\pages 345--348
\endref

\ref \key{14}
\by D.\,Repov\v{s} and P.\,Semenov
\paper On functions of nonconvexity for graphs of continuous functions
\jour J. Math. Anal. Appl.
\vol {\bf 196}
\yr 1995
\pages 1021--1029
\endref

\ref \key{15}
\by D.\,Repov\v{s} and P.\,Semenov
\paper On nonconvexity of graphs of polynomials of several real variables
\jour Set-Valued. Anal.
\vol {\bf 6}
\yr 1998
\pages 39--60
\endref

\ref \key{16}
\by D.\,Repov\v{s} and P.\,Semenov
\paper Selections as an uniform limits of $\d$-continuous $\e$-selections
\jour Set-Valued Anal.
\vol {\bf 7}
\yr 1999
\pages 239--254
\endref


\ref \key{17}
\by D.\,Repov\v{s} and P.\,Semenov
\paper On the relation between the nonconvexity of a set and the nonconvexity of its $\e-$neghborhoods
\jour Math. Notes
\vol {\bf 70}
\yr 2001
\pages 246--259
\endref

\ref \key{18}
\by P. Saveliev
\paper  Fixed points and selections of set-valued maps on spaces with convexity
\jour Int. J. Math. and Math. Sci.
\vol {\bf 24}
\yr 2000
\pages 595--612
\endref

\ref \key{19}
\by  P. V. Semenov
\paper Nonconvexity in problems of multivalued calculus
\jour J. Math. Sci. (N. Y.)
\vol 100
\yr 2000
\pages 2682--2699
\endref

\ref \key{20}
\by H. Torunczyk
\paper Absolute retracts as factors of normed linear spaces
\jour Fund. Math.
\vol {\bf 86 }
\yr 1974
\pages 53--67
\endref

\ref \key{21}
\by M. L. J. van de Vel
\paper A selection theorem for topological convex structures
\jour Trans. Amer. Math. Soc.
\vol {\bf 336}
\yr 1993
\pages 463--496
\endref

\endRefs
\enddocument